\documentclass{article}

\usepackage{amsmath}
\usepackage{amssymb}

\newtheorem{thm}{Theorem}
\newtheorem{conj}{Conjecture}
\newtheorem{lem}{Lemma}

\newtheorem{cor}{Corollary}

\newtheorem{defn}{Definition}

\title{Finding almost squares VI}
\author{Tsz Ho Chan}
\date{}

\begin{document}
\maketitle

\begin{abstract}
In this paper, we continue the study of almost squares and extend the result of the author's fourth paper of the series to almost squares with closer factors.
\end{abstract}

\section{Introduction and Main results}

An almost square is an integer $n$ that can be factored as $n = a b$ with $a, b$ close to $\sqrt{n}$. For example $n = 999999 = 999 \times 1001$ is an almost square. More precisely, for $0 \le \theta \le 1/2$ and $C > 0$,
\begin{defn}
An integer $n$ is a $(\theta, C)$-almost square if $n = a b$ for some integers $a$, $b$ in the interval $[n^{1/2} - C n^\theta, n^{1/2} + C n^\theta]$.
\end{defn}
In [\ref{C}], the author raised the following
\begin{conj} \label{conj1}
Given $1/4 < \theta \le 1/2$, $C > 0$ and any $\epsilon > 0$. For $x$ sufficiently large, almost all intervals $[x - x^{1/2 - \theta + \epsilon}, x + x^{1/2 - \theta + \epsilon}]$ contains a $(\theta, C)$-almost square.
\end{conj}
In fact, one suspects that the above is true without the word ``almost". In [\ref{C}], the author was only able to answer the above question when $\theta = 1/2$ and was not sure how to consider smaller $\theta$. In this paper, we are going to make progress for smaller $\theta$ and get
\begin{thm} \label{mainthm}
Given $1/4 < \theta \le 1/2$, $\epsilon > 0$ and $C > 0$, and let $X > 0$ be a sufficiently large real number. Then the interval $[x, x + x^{1 - 2 \theta} \log^{5+\epsilon} x]$ contains a $(\theta, C)$-almost square for almost all $x \in [X, 2X]$. Here almost all means apart from a set of measure $o(X)$.
\end{thm}
Note that the exponent here is twice that of Conjecture \ref{conj1}. The key idea of its proof is to use a shorter interval of integration which gives rise to a bigger exponent. We also improve a result in [\ref{C}].
\begin{cor} \label{otherthm}
Let $\epsilon > 0$ and $X > 0$ be a sufficiently large real number. Then, for almost all $x \in [X, 2X]$, the interval $[x, x + \log^{5 + \epsilon} x]$ contains an integer $n = a b$ with $\frac{1}{2} x^{1/2} \le a, b \le 2 x^{1/2}$. Here almost all means apart from a set of measure $o(X)$.
\end{cor}

\bigskip

{\bf Some Notations} The notations $f(x) = O(g(x))$, $f(x) \ll g(x)$ and $g(x)
\gg f(x)$ are all equivalent to $|f(x)| \leq C g(x)$ for some
constant $C > 0$. Meanwhile $f(x) = o(g(x))$ means that $\lim_{x \rightarrow \infty} \frac{f(x)}{g(x)} = 0$.

\section{Main idea}

Let $X > Y > 0$ be sufficiently large real numbers and $y \in [X, X+Y]$. Let $1/2 \le L \le U \le X^{1/2}$ and $V \ge 2$ be parameters that may depend on $X$, $Y$ but not $y$. Define
\[
N(s) := \sum_{U - L \le n \le U+L} \frac{1}{n^s} \; \; \text{ and } \; \; \Phi(y) := \mathop{\sum_{y \le n n' \le y(1 + 1/V)}}_{U-L \le n \le U+L} 1.
\]
By Perron's formula
\[
\Phi(y) = \frac{1}{2\pi i} \int_{c - iT}^{c + iT} \zeta(s) N(s) \Bigl[ \Bigl(1 + \frac{1}{V}\Bigr)^s - 1 \Bigr] y^s \frac{ds}{s} + O(|R_y| + |R_{(1 + 1/V)y}|)
\]
where $c = 1 + \frac{1}{\log X}$ and
\[
R_x \ll \mathop{\sum_{x/2 < n < 2x}}_{n \neq x} a_n \min \Bigl(1, \frac{x}{T |x - n|} \Bigr) + \frac{(4x)^c}{T} \sum_{n = 1}^{\infty} \frac{a_n}{n^c}
\]
where $a_n = \sum_{m | n, U - L \le m \le U + L} 1$. Now we shift the line of integration to the left. By Cauchy's reside theorem,
\[
\frac{y}{V} N(1) = \frac{1}{2\pi i} \Bigl( \int_{c-iT}^{c+iT} + \int_{c+iT}^{\eta+iT} + \int_{\eta+iT}^{\eta-iT} + \int_{\eta-iT}^{c-iT} \Bigr) \zeta(s) N(s) \Bigl[ \Bigl(1 + \frac{1}{V} \Bigr)^s - 1 \Bigr] y^s \frac{ds}{s}
\]
for some $1/2 \le \eta < 1$. Thus
\[
\Phi(y) - \frac{y}{V} N(1) = \frac{1}{2\pi i} \Bigl( \int_{c+iT}^{\eta+iT} + \int_{\eta+iT}^{\eta-iT} + \int_{\eta-iT}^{c-iT} \Bigr) \zeta(s) N(s) \Bigl[ \Bigl(1 + \frac{1}{V} \Bigr)^s - 1 \Bigr] y^s \frac{ds}{s} + O(|R_y| + |R_{(1 + 1/V)y}|).
\]
Since $\zeta(\sigma + it) \ll (|t| + 2)^{(1 - \sigma)/3} \log |t|$ for $0 \le \sigma \le 1$ and $|t| \ge 2$ (see [\ref{I}, Theorem 1.9] for example), we have
\[
\int_{c + iT}^{\eta + iT} \zeta(s) N(s) \Bigl[ \Bigl(1 + \frac{1}{V} \Bigr)^s - 1 \Bigr] y^s \frac{ds}{s} \ll \frac{\log T}{T} \int_{\eta}^{c} T^{(1-\sigma)/3} \frac{L}{U^\sigma} y^\sigma d\sigma \ll \frac{L y}{U T} \log T
\]
provided $2 \le \frac{y}{U T^{1/3}} \le X$ and $T \le X$. Therefore as $1 \le L \le U$,
\begin{align} \label{start}
\Phi(y) - \frac{y}{V} N(1) =& \frac{1}{2\pi i} \int_{\eta+iT}^{\eta-iT} \zeta(s) N(s) \Bigl[ \Bigl(1 + \frac{1}{V} \Bigr)^s - 1 \Bigr] y^s \frac{ds}{s} \notag \\
&+ O\Bigl(\frac{L X}{U T} \log T + |R_y| + |R_{(1 + 1/V)y}| \Bigr).
\end{align}
Our goal is to prove that
\[
I_{X,Y} := \frac{1}{Y} \int_{X}^{X+Y} \Big| \Phi(y) - \frac{y}{V} N(1) \Big|^2 dy
\]
is small which would imply that $\Phi(y) \neq 0$ for almost all $y \in [X, X+Y]$.

\section{Some Lemmas}

\begin{lem} \label{diagonal}
For $0 < L \le U/2$,
\[
S_1 = \sum_{N_1 \le n_1 < 2N_1} \sum_{N_2 \le n_2 < 2N_2} \mathop{\sum_{U-L \le m_1, m_2 \le U+L}}_{n_1 m_1 = n_2 m_2} \frac{1}{n_1^{1/2} n_2^{1/2} m_1^{1/2} m_2^{1/2}} \ll \frac{L + U^{1/2}}{U} \log^2 (N_1 N_2 U).
\]
\end{lem}

Proof: Without loss of generality, we may assume that $N_1 / 8 \le N_2 \le 8 N_1$ for otherwise say $N_2 < N_1 / 8$, then since we want $n_1 m_1 = n_2 m_2$,
\[
4 < \frac{N_1}{2 N_2} < \frac{n_1}{n_2} = \frac{m_2}{m_1} \le \frac{U+L}{U-L} \le \frac{3U/2}{U/2} = 3
\]
which is impossible (the other case $N_2 > 8 N_1$ is similar). Thus
\begin{align*}
S_1 \ll& \frac{1}{N_1 U} \sum_{N_1 \le n_1 < 2N_1} \sum_{N_2 \le n_2 < 2N_2} \mathop{\sum_{U-L \le m_1, m_2 \le U+L}}_{n_1 m_1 = n_2 m_2} 1 \\
\le& \frac{1}{N_1 U} \sum_{N_1 \le n_1 < 2N_1} \sum_{U-L \le m_1 \le U+L} d(n_1 m_1) \\
\ll& \frac{1}{N_1 U} \sum_{N_1 \le n_1 < 2N_1} d(n_1) \sum_{U-L \le m_1 \le U+L} d(m_1) \\
\ll& \frac{1}{N_1 U} (N_1 \log N_1) (L \log U + U^{1/2}) \ll \frac{L + U^{1/2}}{U} \log^2 N_1 U
\end{align*}
which gives the lemma.

\begin{lem} \label{nondiagonal}
For $U^\beta < L \le U/2$ with some $0 < \beta < 1/2$,
\begin{align*}
S_2 =& \sum_{N_1 \le n_1 < 2N_1} \sum_{N_2 \le n_2 < 2N_2} \mathop{\sum_{U-L \le m_1, m_2 \le U+L}}_{n_1 m_1 \neq n_2 m_2} \frac{1}{n_1^{1/2} n_2^{1/2} m_1^{1/2} m_2^{1/2}} \frac{1}{|\log (\frac{n_2 m_2}{n_1 m_1})|} \\
\ll& \frac{N_1^{1/2} N_2^{1/2} L^2}{U} \log N_1 N_2 U + \frac{N_1 N_2 L^2}{U^2} \log N_1 N_2 U.
\end{align*}
\end{lem}

Proof: Suppose $N_2 < N_1/8$, then $n_2 m_2 < 3 N_1 U / 8 \le 3 n_1 m_1 / 4$. So $\frac{n_2 m_2}{n_1 m_1} < \frac{3}{4}$ and $|\log \frac{n_2 m_2}{n_1 m_1}| > \log \frac{4}{3}$. Hence
\[
S_2 \ll \sum_{N_1 \le n_1 < 2N_1} \sum_{N_2 \le n_2 < 2N_2} \mathop{\sum_{U-L \le m_1, m_2 \le U+L}}_{n_1 m_1 \neq n_2 m_2} \frac{1}{n_1^{1/2} n_2^{1/2} m_1^{1/2} m_2^{1/2}} \ll \frac{N_1^{1/2} N_2^{1/2} L^2}{U}.
\]
Similarly, we get the same upper bound when $N_2 > 8 N_1$. Thus we may assume $N_1 / 8 \le N_2 \le 8 N_1$. By symmetry, we may even assume $N_1 / 8 \le N_2 \le N_1$. We break down the sum according to
\[
n_1 m_1 - n_2 m_2 = h \; \; \text{ with } \; \; -4 N_2 U \le h \le 4 N_1 U \text{ and } h \neq 0.
\]
Let $d = (m_1, m_2) | h$. Then
\[
x_0 \frac{m_1}{d} - y_0 \frac{m_2}{d} = 1 \text{ for some integers } x_0 \text{ and } y_0
\]
with $- \frac{m_2}{d} \le x_0 \le -1$ and $-\frac{m_1}{d} \le y_0 \le -1$ if $h < 0$; or $1 \le x_0 \le \frac{m_2}{d}$ and $1 \le y_0 \le \frac{m_1}{d}$ if $h > 0$. Let $m_1' = m_1/d$, $m_2' = m_2/d$ and $h' = h/d$. Then
\[
x_0 m_1' - y_0 m_2' = 1
\]
and
\[
n_1 = m_2' t + h' x_0, \; \; n_2 = m_1' t + h' y_0 \text{ for some positive integer } t.
\]
Since $N_1 \le n_1 < 2N_1$ and $N_2 \le n_2 < 2N_2$, we have
\[
\max\Bigl( \frac{N_1}{m_2'} - \frac{h' x_0}{m_2'}, \frac{N_2}{m_1'} - \frac{h' y_0}{m_1'}\Bigr) \le t \le \min\Bigl( \frac{2 N_1}{m_2'} - \frac{h' x_0}{m_2'}, \frac{2 N_2}{m_1'} - \frac{h' y_0}{m_1'}\Bigr).
\]
Hence
\begin{equation} \label{first}
S_2 \ll \frac{1}{N_1^{1/2} N_2^{1/2} U} \sum_{d} \sum_{h'} \mathop{\sum_{(U-L)/d \le m_1', m_2' \le (U+L)/d}}_{(m_1', m_2') = 1} \sum_{t} \frac{1}{\big| \log \bigl(\frac{(m_1' t + h' y_0) m_2'}{(m_2' t + h' x_0) m_1'}\bigr) \big| }
\end{equation}
where the sum over $d$ is from $1$ to $2 U$ (for otherwise $(U+L)/d < 1$), the sum over $h' \neq 0$ is from $-4 N_2 U /d$ to $4 N_1 U/d$, and the sum over $t$ is subjected to the condition above. First, let us separate the contribution from those $d > 4L$. The interval $[(U-L)/d, (U+L)/d]$ has length $2L/d < 1$. Hence $m_1' = m_2' = 1$. So, in this case, we must have $x_0 = 1$ and $y_0 = 0$, and $(U-L)/d \le 1 \le (U+L)/d$ meaning that $U-L \le d \le U+L$. Then the contribution from these $d$'s is
\begin{equation} \label{larged}
\ll \frac{1}{N_1^{1/2} N_2^{1/2} U} \sum_{d} \sum_{h'} \sum_{N_1 - h' \le t \le 2N_1 - h'} \frac{1}{\big| \log \bigl(\frac{t}{t + h'}\bigr) \big| }.
\end{equation}
Suppose $0 < h' < 4 N_1 U / d$ (the other case is similar).
\[
\frac{t}{t + h'} = 1 - \frac{h'}{t + h'} \text{ and } 0 \le \frac{h'}{t + h'} < 1 \text{ when } h' \le \frac{N_1 U}{4 d}.
\]
By $\log (1 - x) \le -x$ when $0 \le x < 1$, we have
\[
\frac{1}{\big| \log \bigl(\frac{t}{t + h'}\bigr) \big| } \le \frac{1}{\frac{h'}{t + h'}} \ll \frac{N_1 U}{h' d}
\]
when $0 < h' \le N_1 U /(4 d)$. When $N_1 U / (4 d) < h' < 4 N_1 U / d$,
\[
\frac{1}{\big| \log \bigl(\frac{t}{t + h'}\bigr) \big| } = \frac{1}{\log (1 + \frac{h'}{t})} \le \frac{1}{\log (1 + \frac{N_1 U/(4d)}{4 U N_2/d})} \ll 1 \ll \frac{N_1 U}{h' d}.
\]
Thus
\begin{equation} \label{larged2}
(\ref{larged}) \ll \frac{1}{U} \sum_{U-L \le d \le U+L} \sum_{0 < h' \le 4 N_1 U / d} \frac{N_1 U}{d h'} \ll \frac{N_1 L}{U} \log N_1 U.
\end{equation}

From now on, we restrict our attention to (\ref{first}) with $d \le 4L$. Suppose $-4 N_2 U \le h' < 0$ (the other case is similar). Observe that
\[
\frac{(m_1' t + h' y_0) m_2'}{(m_2' t + h' x_0) m_1'} = 1 + \frac{-h'}{(m_2' t + h' x_0) m_1'} \text{ and } \frac{N_1 U}{2 d} \le (m_2' t + h' x_0) m_1' = n_1 m_1' \le \frac{4 N_1 U}{d}.
\]
If $-\frac{N_2 U}{2 d} \le h' < 0$, then $0 < \frac{-h'}{(m_2' t + h' x_0) m_1} \le 1$. By $\log (1 + x) \ge \frac{x}{2}$ for $0 \le x \le 1$,
\[
\frac{1}{\big| \log \bigl(\frac{(m_1' t + h' y_0) m_2'}{(m_2' t + h' x_0) m_1'}\bigr) \big| } \le
\frac{2}{\frac{-h'}{(m_2' t + h' x_0) m_1}} \ll \frac{N_1 U}{|h'| d}.
\]
If $- \frac{4 N_2 U}{d} \le h' < -\frac{N_2 U}{2 d}$, then $1 < \frac{-h'}{(m_2' t + h' x_0) m_1'} \le 8$. Hence
\[
\frac{1}{\big| \log \bigl(\frac{(m_1' t + h' y_0) m_2'}{(m_2' t + h' x_0) m_1'}\bigr) \big| } \le \frac{1}{\log 2} \ll \frac{N_1 U}{|h'| d}.
\]
Combining the above estimates with (\ref{first}), (\ref{larged}) and (\ref{larged2}), we have
\[
S_2 \ll  \frac{N_1 L}{U} \log N_1 U + \frac{1}{N_1^{1/2} N_2^{1/2} U} \sum_{U/4N_1 < d \le 4L} \sum_{h'} \mathop{\sum_{(U-L)/d \le m_1', m_2' \le (U+L)/d}}_{(m_1', m_2') = 1} \frac{d N_1}{U} \frac{N_1 U}{|h'| d}.
\]
The reason we have $U/4N_1 < d$ is that if $d \le U/4N_1$, then $t < 2N_1 / m_2' \le 2N_1 / (U/2d) = 4N_1 d / U \le 1$ as $L \le U/2$. Hence
\begin{align*}
S_2 \ll& \frac{N_1 L}{U} \log N_1 U + \frac{N_1^{1/2} N_2^{1/2}}{U} \sum_{U/4N_1 < d \le 4L} \sum_{h'} \frac{1}{|h'|} \mathop{\sum_{(U-L)/d \le m_1', m_2' \le (U+L)/d}}_{(m_1', m_2') = 1} 1 \\
\ll& \frac{N_1^{1/2} N_2^{1/2} L}{U} \log N_1 N_2 U + \frac{N_1 N_2 L^2}{U^2} \log N_1 N_2 U
\end{align*}
which gives the lemma.

\begin{lem} \label{keylemma}
For $U^\beta < L \le U/2$ with $1/2 \le \beta \le 1$,
\[
I := \int_{1}^{T} \Big| \zeta \Bigl(\frac{1}{2} + it \Bigr) N \Bigl(\frac{1}{2} + it \Bigr) \Big|^2 dt \ll \frac{T L}{U} \log^2 TU + \frac{T^{1/2} L^2}{U} \log TU.
\]
\end{lem}

Proof: First, recall the approximate functional equation of $\zeta(s)$ (see [\ref{I}, Theorem 4.1] for example):
\[
\zeta(\sigma + it) = \sum_{n_1 \le \sqrt{t/2\pi}} \frac{1}{n_1^{\sigma + it}} + \chi(\sigma + it) \sum_{n_2 \le \sqrt{t/2\pi}} \frac{1}{n_2^{1-\sigma-it}}
+ O(t^{-1/4})
\]
where
\[
\chi(s) = \frac{(2\pi)^s}{2 \Gamma(s) \cos(\pi s/2)}.
\]
Put $\sigma = 1/2$ and use $|a + b + c|^2 \le 3(|a|^2 + |b|^2 + |c|^2)$, we get
\begin{align*}
I \ll& \int_{1}^{T} \Big| \sum_{n_1 \le \sqrt{t/2\pi}} \frac{1}{n_1^{1/2 + it}} N \Bigl(\frac{1}{2} + it\Bigr) \Big|^2 dt \\
&+ \int_{1}^{T} \Big|\chi\Bigl(\frac{1}{2} + it \Bigr) \Big|^2 \Big| \sum_{n_2 \le \sqrt{t/2\pi}} \frac{1}{n_2^{1/2 - it}} N \Bigl( \frac{1}{2} + it \Bigr)
\Big|^2 dt \\
&+ \int_{1}^{T} \frac{1}{t^{1/2}} \Big| N \Bigl(\frac{1}{2} + it \Bigr) \Big|^2 dt =: I_1 + I_2 + I_3.
\end{align*}
We estimate $I_3$ first. By integration by parts,
\[
I_3 = \frac{I(T)}{T^{1/2}} + \frac{1}{2} \int_{1}^{T} \frac{I(u)}{u^{3/2}} du
\]
where $I(u) = \int_{1}^{u} |N(\frac{1}{2} + it)|^2 dt$. By Montgomery and Vaughan's mean value theorem [\ref{MV}],
\[
I(u) \ll u \frac{L}{U} + L.
\]
Hence
\[
I_3 \ll T^{1/2} \frac{L}{U} + L.
\]
Next, we estimate $I_1$. By expanding things out, we have
\begin{align*}
I_1 =& \int_{1}^{T} \sum_{n_1, n_2 \le \sqrt{t/2\pi}} \sum_{U-L \le m_1, m_2 \le U+L} \frac{1}{n_1^{1/2} n_2^{1/2} m_1^{1/2} m_2^{1/2}} \Bigl(\frac{n_2 m_2}{n_1 m_1} \Bigr)^{it} dt \\
=& \sum_{n_1, n_2 \le \sqrt{T/2\pi}} \sum_{U-L \le m_1, m_2 \le U+L} \frac{1}{n_1^{1/2} n_2^{1/2} m_1^{1/2} m_2^{1/2}}
\int_{2 \pi \max(n_1^2, n_2^2)}^{T} \Bigl(\frac{n_2 m_2}{n_1 m_1} \Bigr)^{it} dt \\
\ll& T \mathop{\sum_{n_1, n_2 \le \sqrt{T/2\pi}} \sum_{U-L \le m_1, m_2 \le U+L}}_{n_1 m_1 = n_2 m_2} \frac{1}{n_1^{1/2} n_2^{1/2} m_1^{1/2} m_2^{1/2}} \\
&+ \mathop{\sum_{n_1, n_2 \le \sqrt{T/2\pi}} \sum_{U-L \le m_1, m_2 \le U+L}}_{n_1 m_1 \neq n_2 m_2} \frac{1}{n_1^{1/2} n_2^{1/2} m_1^{1/2} m_2^{1/2}}
\Big| \frac{1}{\log \bigl(\frac{n_2 m_2}{n_1 m_1}\bigr)} \Big|.
\end{align*}
Apply Lemma \ref{diagonal} and Lemma \ref{nondiagonal}, we obtain
\[
I_1 \ll T \frac{L + U^{1/2}}{U} \log^2 TU + \frac{T^{1/2} L^2}{U} \log TU + \frac{T L^2}{U^2} \log TU.
\]
As $|\chi(1/2 + it)| = 1$,
\[
I_2 \ll T \frac{L + U^{1/2}}{U} \log^2 TU + \frac{T^{1/2} L^2}{U} \log TU + \frac{T L^2}{U^2} \log TU
\]
by almost the same argument as $I_1$. Combining the bounds for $I_1$, $I_2$ and $I_3$, we have the lemma.

\bigskip

Finally, we need a lemma to bound the error terms in Perron's formula.
\begin{lem} \label{perron}
For $\epsilon > 0$, $1 \le a \le 2$ and $1 \le T,  Y \le X$,
\[
\frac{1}{Y} \int_{X}^{X+Y} |R_{ay}|^2 dy \ll_\epsilon \frac{L^2 X^2}{U^2 T^2} \log^2 X + \frac{X L^2}{T Y} \log^2 X.
\]
\end{lem}

Proof: Recall
\[
R_x \ll \mathop{\sum_{x/2 < n < 2x}}_{n \neq x} a_n \min \Bigl(1, \frac{x}{T |x-n|} \Bigr) + \frac{(4x)^c}{T} \sum_{n = 1}^{\infty} \frac{a_n}{n^c}
\]
where $c = 1 + 1/\log X$ and $a_n = \sum_{m | n, U - L \le m \le U + L} 1$. Splitting the sum and using $\zeta(s) = \frac{1}{s-1} + O(1)$,
\begin{align*}
R_{ay} \ll& \sum_{ay - \frac{ay}{T} \le n \le ay + \frac{ay}{T}} \mathop{\sum_{m | n}}_{U - L \le m \le U + L} 1 + \frac{y}{T} \sum_{\frac{ay}{2} < n < ay - \frac{ay}{T}} \frac{1}{ay - n} \mathop{\sum_{m | n}}_{U - L \le m \le U + L} 1 \\
&+ \frac{y}{T} \sum_{\frac{ay}{2} < n < ay - \frac{ay}{T}} \frac{1}{n - ay} \mathop{\sum_{m | n}}_{U - L \le m \le U + L} 1 + \frac{X}{T} \sum_{n = 1}^{\infty} \frac{a_n}{n^c} \mathop{\sum_{m | n}}_{U - L \le m \le U + L} 1 \\
\ll& \sum_{U - L \le m \le U + L} \sum_{\frac{ay}{m} - \frac{ay}{Tm} \le n' \le \frac{ay}{m} + \frac{ay}{Tm}} 1 + \frac{y}{T} \sum_{U - L \le m \le U + L} \sum_{\frac{ay}{2m} < n' < \frac{ay}{m} - \frac{ay}{Tm}} \frac{1}{ay - m n'} \\
&+ \frac{y}{T} \sum_{U - L \le m \le U + L} \sum_{\frac{ay}{2m} < n' < \frac{ay}{m} - \frac{ay}{Tm}} \frac{1}{m n' - a y} + \frac{X}{T} \sum_{U - L \le m \le U + L} \frac{1}{m} \sum_{n' = 1}^{\infty} \frac{1}{n'^c} \\
\ll& \sum_{U - L \le m \le U + L} \sum_{\frac{ay}{m} - \frac{ay}{Tm} \le n' \le \frac{ay}{m} + \frac{ay}{Tm}} 1 + \frac{L X}{U T} \log X.
\end{align*}
Let $f_\Delta (x) = \max (1 - |x|/\Delta, 0)$ where $\Delta = \frac{10X}{UT}$ and $g_\Delta (x) = \sum_{n = -\infty}^{\infty} f_\Delta (x - n)$. Hence
\[
R_{ay} \ll \sum_{U - L \le m \le U + L} g_\Delta \Bigl(\frac{ay}{m}\Bigr) + \frac{L X}{U T} \log X.
\]
Now we make use of the Fourier series of $g_\Delta(x)$, say
\begin{align*}
g_\Delta (x) =& \sum_{k = -\infty}^{\infty} a_k e(k x) \text{ where } a_k = \Delta \Bigl(\frac{\sin \pi k \Delta}{\pi k \Delta} \Bigr)^2 \\
=& \Delta + \sum_{0 < |k| \le K} a_k e(k x) + O (\Delta)
\end{align*}
where $K = \frac{1}{\Delta^2}$. Therefore
\begin{equation} \label{i1}
\int_{X}^{X+Y} |R_{ay}|^2 dy \ll Y \frac{L^2 X^2}{U^2 T^2} \log^2 X + \int_{X}^{X+Y} \Big| \sum_{0 < k \le K} a_k \sum_{U - L \le m \le U + L} e \Bigl( \frac{k a y}{m} \Bigr) \Big|^2 dy.
\end{equation}
It remains to deal with the last integral denoted by $J$. Expanding things out, isolating the diagonal terms, and interchanging summations and integration, we have
\begin{align*}
J \ll& Y \sum_{0 < k, l \le K} a_k a_l \mathop{\sum_{U - L \le m, n \le U + L}}_{k n = l m} 1 + U^2 \sum_{0 < k, l \le K} a_k a_l \mathop{\sum_{U - L \le m, n \le U + L}}_{k n \neq l m} \frac{1}{|k n - l m|} \\
\ll& Y \sum_{0 < k, l \le K} a_k a_l \mathop{\sum_{U - L \le m, n \le U + L}}_{k n = l m} 1 + U^2 \sum_{0 < t \le K (U+L)} \frac{1}{t} \mathop{\sum_{U - L \le m, n \le U + L}}_{k n - l m = t} 1 =: S_1 + S_2.
\end{align*}
First, we deal with $S_1$. Suppose $d = (m, n)$. Let $m' = m/d$ and $n' = n/d$. Then $k = m' s$ and $l = n' s$ for positive integer $s$. Then
\begin{align*}
S_1 =& Y \sum_{d \le 2L} \sum_{U/d - L/d \le m', n' \le U/d + L/d} \sum_{s = 1}^{\infty} a_{m' s} a_{n' s} \\
\ll& Y \sum_{d \le 2L} \sum_{U/d - L/d \le m', n' \le U/d + L/d} \sum_{s \le \frac{d}{\Delta U}} \Delta^2 + Y \sum_{d \le 2L} \sum_{U/d - L/d \le m', n' \le U/d + L/d} \sum_{s > \frac{d}{\Delta U}} \frac{d^4}{\Delta^2 U^4} \frac{1}{s^4} \\
\ll& Y \frac{\Delta L^2}{U} \log X.
\end{align*}
For $S_2$, suppose $d = (m,n)$ which divides $t$. Let $m' = m/d$, $n' = n/d$ and $t' = t/d$. For each pair of $m'$ and $n'$, say $x_0 n' - y_0 m' = 1$ for some $0 \le x_0 \le m'$ and $0 \le y_0 \le n'$. Then $k = m' s + x_0 t'$ and $l = n' s + y_0 t'$ for non-negative integer $s$. Thus
\begin{align*}
J \ll& Y \frac{\Delta L^2}{U} \log X + U^2 \sum_{0 < t' \le K (U+L)} \frac{1}{t'} \sum_{d \le 2L} \frac{1}{d} \sum_{U/d - L/d \le m' \le n' \le U/d + L/d} \sum_{s = 1}^{\infty} a_{m' s + x_0 t'} a_{n' s + y_0 t'} \\
\ll& Y \frac{\Delta L^2}{U} \log X + U^2 \sum_{0 < t' \le K (U+L)} \frac{1}{t'} \sum_{d \le 2L} \frac{1}{d} \sum_{U/d - L/d \le m' \le n' \le U/d + L/d} \sum_{s \le 1/(\Delta m')} \Delta^2 \\
&+ U^2 \sum_{0 < t' \le K (U+L)} \frac{1}{t'} \sum_{d \le 2L} \frac{1}{d} \sum_{U/d - L/d \le m' \le n' \le U/d + L/d} \sum_{s > 1/(\Delta n')} \frac{1}{\Delta (m' s)^2} \frac{1}{\Delta (n' s)^2} \\
\ll& Y \frac{\Delta L^2}{U} \log X + \Delta U L^2 \log^2 X.
\end{align*}
Putting this into (\ref{i1}) and dividing by $Y$, we get the lemma.

\section{Proof of Theorem \ref{mainthm} and \ref{otherthm}}

\[
J_{X,Y} := \frac{1}{Y} \int_{X}^{X+Y} \Big| \int_{\eta+iT}^{\eta-iT} \zeta(s) N(s) \Bigl[ \Bigl(1 + \frac{1}{V} \Bigr)^s - 1 \Bigr] y^s \frac{ds}{s} \Big|^2 dy.
\]
Multiplying things out and integrating over $y$, we have
\[
J_{X,Y} \ll \frac{X^2}{Y V^2} \int_{-T}^{T} \int_{-T}^{T} |\zeta(\eta + iu)| |\zeta(\eta + iv)| |N(\eta + iu)| |N(\eta + iv)| \frac{du dv}{1 + |u - v|}.
\]
Here we use the fact that $(1 + \frac{1}{V})^{\eta + iu} = e^{(\eta + iu) \log(1 + 1/V)} = e^{O(|\eta + iu|/V)} = 1 + O(\frac{|\eta + iu|}{V})$ when $|u| \le V/2$ and $|\frac{(1 + 1/V)^{\eta + iu} - 1}{\eta + iu}| \le \frac{(1 + 1/V)^\eta + 1}{V/2} \ll \frac{1}{V}$ when $|u| > V/2$. As $2abcd \le (ac)^2 + (bd)^2$,
\begin{align*}
J_{X,Y} \ll& \frac{X^2}{Y V^2} \int_{-T}^{T} |\zeta(\eta + iu)|^2 |N(\eta + iu)|^2 \int_{-T}^{T} \frac{1}{1 + |u - v|} dv du \\
\ll& \frac{X^2}{Y V^2} \log T \int_{-T}^{T} |\zeta(\eta + iu)|^2 |N(\eta + iu)|^2 du \\
\ll& \frac{X^2}{Y V^2} \log T \Bigl( \frac{T L}{U} \log^2 TU + \frac{T^{1/2} L^2}{U} \log TU \Bigr)
\end{align*}
by Lemma \ref{keylemma}. Hence, together with Lemma \ref{perron}, we have
\[
I_{X,Y} \ll \frac{X^2}{Y V^2} \log X \Bigl(\frac{T L}{U} \log^2 X + \frac{T^{1/2} L^2}{U} \log X + \frac{T L^2}{U^2} \log X \Bigr) + \frac{L^2 X^2}{U^2 T^2} \log^2 X + \frac{X L^2}{T Y} \log^2 X
\]
as $U \le X^{1/2}$ and $T \le X$. Let
\[
\mathcal{B} := \{y \in [X, X+Y]: \Phi(y) = 0\} \text{ and } |\mathcal{B}| \text{ be its measure.}
\]
Then
\[
\frac{|\mathcal{B}|}{Y} \frac{X^2}{V^2} \Bigl( \frac{L}{U} \Bigr)^2 \ll \frac{X^2}{Y V^2} \log X \Bigl(\frac{T L}{U} \log^2 X + \frac{T^{1/2} L^2}{U} \log X \Bigr) + \frac{L^2 X^2}{U^2 T^2} \log^2 X + \frac{X L^2}{T Y} \log^2 X.
\]
Let $U = X^{1/2}$ and $L = \frac{C X^\theta}{2C + 3}$ with $1/4 < \theta \le 1/2$,
\begin{equation} \label{measure}
|\mathcal{B}| \ll T \frac{U}{L} \log^3 X + T^{1/2} U \log^2 X + \frac{Y V^2}{T^2} \log^2 X + \frac{V^2 U^2}{X T} \log^2 X.
\end{equation}
Set $Y = X^{1/2} L \le C X^{1/2 + \theta} / 3$. By picking $T = X^{2 \theta} / \log^{4 + \epsilon/2} X$ and $V = X^{2 \theta} / \log^{5 + \epsilon} X$, one can check that $|\mathcal{B}| = o(Y)$ and $X^{1/2} - C X^\theta < n' < X^{1/2} + C X^\theta$ as $\theta > 1/4$ and $\epsilon$ can be arbitrarily small. This proves Theorem \ref{mainthm}.


White Station High School \\
514 S. Perkins Road, \\
Memphis, TN 38117 \\
U.S.A.

\end{document}